\documentclass[11pt,bezier,amstex]{article}  
\usepackage[textsize=small]{todonotes}       

\usepackage{color}              
\usepackage{graphicx,setspace}
\usepackage[]{amsmath}
\usepackage{amssymb}
\usepackage{fullpage}
\usepackage{hyperref}
\usepackage{enumerate, amsfonts, amsthm, bbm}

\definecolor{MyDarkBlue}{rgb}{0,0.08,0.50}
\definecolor{BrickRed}{rgb}{0.65,0.08,0}

\hypersetup{
colorlinks=true,       
    linkcolor=MyDarkBlue,          
    citecolor=BrickRed,        
    filecolor=red,      
    urlcolor=cyan           
}

\numberwithin{equation}{section}

\newcommand{\ra}{\rightarrow}

\usepackage{tabularx}

\newcommand{\pas}[1]{T_x(#1)}
\newcommand{\dd}{{\rm d}}
\newcommand{\tm}{\theta_{{\rm max}}}
\newcommand{\Bs}{B_{\ast}}
\newcommand{\expect}[1]{{\mathbb E}\left[#1\right]}

\newcommand{\prob}[1]{{\mathbb{P}}\left( #1 \right)}

\newcommand{\eqn}[1]{\begin{equation} #1 \end{equation}}
\newcommand{\eqan}[1]{\begin{align} #1 \end{align}}
\newcommand{\e}{{\rm e}}
\newcommand{\nn}{\nonumber}
\newcommand{\lst}[1]{{\varphi}_{#1}}

\newtheorem{theorem}{Theorem}[]

\newtheorem{proposition}[theorem]{Proposition}
\theoremstyle{definition}

\begin{document}
\title{{\Large {\bf First passage times to congested states of many-server systems in the Halfin-Whitt regime}}}

\author{Brian Fralix\footnotemark[1]  \and
Charles Knessl\footnotemark[2] \and
        Johan S.H. van Leeuwaarden\footnotemark[3]
        }

\date{\today}

\maketitle

\footnotetext[1]{Clemson University, Department of Mathematical Sciences, O-110 Martin Hall, Box 340975, Clemson, SC 29634, USA. Email: {\tt bfralix@clemson.edu}}

\footnotetext[2]{University of Illinois at Chicago, Department of Mathematics, Statistics and Computer Science, 815 South Morgan Street, Chicago, IL 60607-7045, USA. Email address: {\tt knessl@uic.edu}}

\footnotetext[3]{Eindhoven University of Technology and EURANDOM, P.O. Box 513, 5600 MB Eindhoven, The Netherlands. Email address: {\tt j.s.h.v.leeuwaarden@tue.nl}}

\begin{abstract}
We consider the heavy-traffic approximation to the $GI/M/s$ queueing system in the Halfin-Whitt regime, where both the number of servers $s$ and the arrival rate $\lambda$ grow large (taking the service rate as unity), with $\lambda=s-\beta\sqrt{s}$ and $\beta$ some constant. In this asymptotic regime, the queue length process can be approximated by a diffusion process that behaves like a Brownian motion with drift above zero and like an Ornstein-Uhlenbeck process below zero. We analyze the first passage times of this hybrid diffusion process to levels in the state space that represent congested states in the original queueing system.

\vspace{1.5mm}


\vspace{1.5mm}

\noindent {\it AMS $2000$ Subject Classification}: 60K25, 60J60, 60J70, 34E05.
\end{abstract}

\section{Introduction}
Halfin and Whitt \cite{halfinwhitt} introduced in their 1981 paper a new heavy-traffic limit theorem for the $GI/M/s$ system. They demonstrated how under certain conditions a sequence of normalized queue-length processes converges to a process that behaves like a Brownian motion with drift above zero and like an Ornstein-Uhlenbeck process below zero. We refer to this hybrid diffusion process as the {\it Halfin-Whitt diffusion}.


In \cite{halfinwhitt} it is established that by setting the traffic intensity $\rho=1-\beta/\sqrt{s}$, $\beta\in(0,\infty)$,
the number of customers in the $M/M/s$ system can be roughly expressed as $s+\sqrt{s}X(t)$ for $s$ sufficiently large and $(X(t))_{t\geq 0}$ the Halfin-Whitt diffusion. The boundary between the Brownian motion and the Ornstein-Uhlenbeck process can be thought of as the number of servers, and $(X(t))_{t\geq 0}$ will keep fluctuating between these two regions. The process mimics a single server queue above zero, and an infinite server queue below zero, for which Brownian motion and the Ornstein-Uhlenbeck process are indeed the respective heavy-traffic limits. As $\beta$ increases, capacity grows and the Halfin-Whitt diffusion will spend more time below zero.

The diffusion process $(X(t))_{t\geq 0}$ can thus be employed to obtain simple approximations for the system behavior. The steady-state characteristics of the diffusion were studied in \cite{halfinwhitt}. It is also of interest to study time-dependent characteristics like the mixing times, time-dependent distributions and first passage times to enhance our understanding of how the $GI/M/s$ system behaves over various time and space scales. The mixing time is closely related to the spectral gap, which for the Halfin-Whitt diffusion $(X(t))_{t\geq 0}$ has been identified by Gamarnik and Goldberg \cite{goldberggamarnik} building on the results of van Doorn \cite{vandoorn} on the spectral gap of the $M/M/s$ queue. An alternative derivation of this spectral gap was presented in \cite{vlk}, along with expressions for the Laplace transform over time, and the
large-time asymptotics for the time-dependent density. In this paper we derive results for first passage times to large levels. Such large levels typically correspond to highly congested states, in which users start receiving degraded service. An expression for the mean first passage time was derived in Maglaras and Zeevi \cite{maglaraszeevi}. We shall derive the Laplace transform of the first passage time density. From this Laplace transform, we can derive not only  all moments, but also expressions for the first passage time density in various asymptotic regimes.

Mathematically, determining the Laplace transform of the first passage or time-dependent distributions for the present diffusion process involves analyzing a Schr\"{o}dinger  type equation with a piecewise parabolic potential function, or, equivalently, a Fokker-Planck equation with a piecewise linear drift. Such problems  arise in a variety of other applications, such as linear systems driven by white noise \cite{At67,AtC}, the Kramers' problem \cite{MP} and escape over potential barriers \cite{LRH}. Invariably, the solution involves the parabolic cylinder functions (see also \cite{a&s,vlka,TEMME} for more background on the parabolic cylinder function). The main results are presented in Section \ref{sec:main} and the proofs are given in Sections \ref{secc3}-\ref{secnew}.

\section{Main results}\label{sec:main}

For the Halfin-Whitt diffusion process, define $\pas{b}$ as the first passage time out of the interval $(-\infty,b)$, starting at $x<b$ with $b>0$. Define
\eqn{
\lst{x,b}(\theta)=\mathbb E[\e^{-\theta \pas{b}}],
}
so that if $P(x,t){\rm d}t=\prob{\pas{b}\in[t,t+{\rm d}t]}$
\eqn{\label{lst}
\lst{x,b}(\theta)=\int_{0}^\infty {\rm e}^{-\theta t}P(x,t){\rm d}t, \quad \Re(\theta)>0.
}

The first passage time density $P$ satisfies the backward Kolmogorov equation
\begin{equation}\label{diffusioneq}
P_t=A(x)P_x+\tfrac{1}{2}B(x)P_{xx}; \quad x<b, \ t>0
\end{equation}
with $P(b,t)=\delta(t)$ (the Dirac function) and
\begin{equation}\label{eqr2.1}
A(x)=\left\{
       \begin{array}{ll}
         -\beta, & \hbox{$x> 0$}, \\
         -x-\beta, & \hbox{$x<0$.}
       \end{array}
     \right.
\end{equation}
We also require  $P$ and $P_x$ to be continuous at $x=0$. Here the diffusion coefficient is $B(x)=1+c^2$ where $c$ is the coefficient of variation for the interarrival distribution of the $GI/M/s$ system. For $GI$=$M$ we have $c=1$, and in general we can rescale $x$ so as to make $B(x)=2$, which we henceforth assume.

Let $D_\nu(z)$ denote the parabolic cylinder function with index $\nu$ and argument $z$, which is defined, for example, by the integrals
\eqan{
D_\nu(z)&=\frac{{\rm e}^{-z^2/4}}{\Gamma(-\nu)}\int_0^\infty {\rm e}^{-zu}{\rm e}^{-u^2/2}u^{-\nu-1}{\rm d}u, \quad \Re(\nu)<0,\label{Aa}\\
D_\nu(z)&=\frac{{\rm e}^{z^2/4}}{i\sqrt{2\pi}}\int_\mathcal{C} u^\nu {\rm e}^{u^2/2}{\rm e}^{-uz}{\rm d}u.\label{Bb}
}
Here, $\Gamma(\cdot)$ is the Gamma function, and the contour $\mathcal{C}$ in the second integral is a vertical Bromwich contour in the half-plane $\Re(u)>0$. It is well known that $D_\nu(z)$ is an entire function of both index $\nu$ and argument $z$, and various properties of $D_\nu(z)$ are given in \cite[Chapter 19]{a&s} and \cite[p.~1092-1095]{gradshteyn}.

Define \eqan{\label{M}
M(\theta;\beta,b)=\cosh\left(\frac{b}{2}\sqrt{\beta^2+4\theta}\right)-
\frac{2D_{-\theta}'(-\beta)}{D_{-\theta}(-\beta)}\frac{\sinh\left(\frac{b}{2}\sqrt{\beta^2+4\theta}\right)}{\sqrt{\beta^2+4\theta}}
}
with $D_{-\theta}'(-\beta)=-\frac{d}{d\beta}D_{-\theta}(-\beta)$.
Below we give expressions for $\lst{x,b}(\theta)$, where we must distinguish between the cases $x>0$ and $x<0$.

\begin{theorem}\label{thmneg}
Let $x<0$. Then, with $M(\theta;\beta,b)$ as defined in \eqref{M},
\eqan{\label{28}
\lst{x,b}(\theta) =\frac{1}{M(\theta;\beta,b)}\frac{D_{-\theta}(-\beta-x)}{D_{-\theta}(-\beta)}\exp\left(-\frac{\beta (b-x)}{2}+\frac{x^2}{4}\right).
}
\end{theorem}

\begin{theorem}\label{thmpos}
Let $x>0$. Then, with $M(\theta;\beta,b)$ as defined in \eqref{M},
\eqan{\label{29}
\lst{x,b}(\theta) =\exp\left(\frac{\beta (x-b)}{2}\right)\left[\frac{\sinh\left(\frac{x}{2}\sqrt{\beta^2+4\theta}\right)}{\sinh\left(\frac{b}{2}\sqrt{\beta^2+4\theta}\right)}
-\frac{1}{M(\theta;\beta,b)}\frac{\sinh\left(\frac{(x-b)}{2}\sqrt{\beta^2+4\theta}\right)}{\sinh\left(\frac{b}{2}\sqrt{\beta^2+4\theta}\right)}\right].
}

\end{theorem}
Using
\eqn{
\expect{\pas{b}}=-\frac{\dd }{\dd \theta}\lst{x,b}(\theta)\Big|_{\theta=0},
}
we obtain after tedious calculations the following result for the mean first passage time, which is in agreement with the result obtained in a different manner by Maglaras and Zeevi \cite{maglaraszeevi}.

\begin{proposition}{\rm \cite[Proposition 3]{maglaraszeevi}}\label{mean}
If $x>0$ then
\eqn{\label{211}
\expect{\pas{b}}=\frac{x-b}{\beta}+\left(\e^{\beta b}-\e^{\beta x}\right)\left[\frac{1}{\beta^2}+\frac{1}{\beta}\int_{0}^\infty \e^{\beta u-u^2/2}{\rm d}u\right].
}
If $x<0$ then
\eqn{\label{212}
\expect{\pas{b}}=\frac{\e^{\beta b}-1-\beta b}{\beta^2}+\frac{\e^{\beta b}-1}{\beta}\int_{0}^\infty \e^{\beta u-u^2/2}{\rm d}u-\int_{0}^\infty \e^{\beta u-u^2/2}\left(\frac{\e^{ux}-1}{u}\right){\rm d}u.
}
\end{proposition}

Though \eqref{211} and \eqref{212} are already fairly simple, we give some asymptotic formulas below that yields further insight on the magnitude of the mean passage time (the derivation is standard and therefore omitted).
\begin{proposition}\label{meanregime}
{\rm (a)} For $b\to\infty$ and $b-x\to\infty$,
\eqn{\label{star}
\expect{\pas{b}}\sim\e^{\beta b}\left[\frac{1}{\beta^2}+\frac{1}{\beta}\int_{0}^\infty \e^{\beta u-u^2/2}{\rm d}u\right].
}
If $b-x=O(1)$ the above term should be multiplied by $1-\e^{-\beta(b-x)}$. If $b,\beta\to\infty$ then \eqref{star} simplifies further to $\expect{\pas{b}}\sim\sqrt{2\pi}\beta^{-1}\e^{\beta b}\e^{\beta^2/2}$.\\
{\rm (b)}
For $\beta\to -\infty$ and $x,b=O(|\beta|)$ {\rm(}possibly $o(|\beta|)${\rm)},
\eqn{\label{starr}
\expect{\pas{b}}\sim\left\{
                      \begin{array}{ll}
                        \frac{b-x}{-\beta}, & \hbox{$x\in[0,b)$,} \\
                        \frac{b}{-\beta}+\log\Big(1+\frac{x}{\beta}\Big), & \hbox{$x\in(-\infty,0]$.}
                      \end{array}
                    \right.
}
{\rm (c)}
For $\beta\to \infty$ with $\beta=O(b^{-1})$ and $x=O(b)$,
\eqn{\label{starrr}
\expect{\pas{b}}\sim\left\{
                      \begin{array}{ll}
                        \beta^{-2}[\e^{\beta b}-\e^{\beta x}+\beta(x-b)], & \hbox{$x/b\in[0,1)$,} \\
                        \beta^{-2}[\e^{\beta b}-1-\beta b], & \hbox{$x/b<0$.}
                      \end{array}
                    \right.
}
\end{proposition}

We note that in \eqref{star} the mean first passage time is exponentially large and independent of the starting point $x$, in \eqref{starr} it is asymptotically $O(1)$, and \eqref{starrr} represents the transition between these two cases, where $\expect{\pas{b}}=O(b^2)$.

The Laplace transform $\lst{x,b}(\theta)$ is analytic in the entire $\theta$-plane, except for singularities in the range $\Re(\theta)<0$. Hence, the asymptotic behavior of $\pas{b}$ is determined by the  singularity $\tm$ closest to the imaginary axis. In fact, from  Theorems \ref{thmneg} and \ref{thmpos}
 it follows that $\tm$ will be the largest negative solution to
\begin{equation}\label{key}
D_{-\theta}(-\beta)M(\theta;\beta,b)=0.
\end{equation}
Note that \eqref{M} and \eqref{29} are invariant under the change $\sqrt{\beta^2+4\theta}\rightarrow-\sqrt{\beta^2+4\theta}$, so there is no branch point at  $\theta=-\beta^2/4$. It seems impossible to find a closed-form solution to \eqref{key}. We therefore consider several asymptotic regimes:
\begin{itemize}
\item[(i)] Large levels: $b\ra\infty$ and $\beta$ fixed.
\item[(ii)] Large levels and over/undercapacity: $\beta\ra\pm \infty,b\ra\infty$ at the same rate ($|\beta|/b$ fixed).
\item[(iii)] Small levels and undercapacity: $\beta\ra - \infty$ and $b\ra 0$.
\end{itemize}

Regime (i) represents the situation of reaching highly congested states, corresponding to large levels $b$.
We have the following results:
\begin{proposition}[Regime (i)]\label{prop4}
If $\beta<0$ is fixed, and $b\rightarrow\infty$, then
\eqan{\label{214}
\tm = -\frac{1}{4}\beta^2-\frac{\pi^2}{b^2}\left[1+\frac{2}{b}\frac{D_{\beta^2/4}(-\beta)}{D_{\beta^2/4}'(-\beta)}+O(b^{-2})\right].
}
If $\beta=0$, and $b\rightarrow\infty$, then
\eqan{\label{215}
\tm = -\frac{\pi^2}{4b^2}\left[1-\frac{\sqrt{2\pi}}{b}+O(b^{-2})\right].
}
If $\beta>0$ is fixed, and $b\rightarrow\infty$, then
\eqan{\label{3rdr}
\tm \sim -\frac{\beta^2\e^{-\beta b}}{1+\beta\e^{\beta^2/2}\int_{-\infty}^\beta \e^{-u^2/2}{\rm d}u}.
}
\end{proposition}

We can generalize \eqref{215} to the case where $b\rightarrow\infty$ and $\beta\ra 0$ with $\beta b=\gamma$ fixed, where we have
\eqan{\label{stst}
\tm\sim -\frac{1}{\beta^2}\left[\frac{\gamma^2}{4}+\omega(\gamma)\right]
}
where $\omega$ is the solution of the smallest absolute value to
\eqan{
\frac{\tan(\sqrt{\omega})}{\sqrt{\omega}}=\frac{\tanh(\sqrt{-\omega})}{\sqrt{-\omega}}=\frac{2}{\gamma}.
}
It follows that if $\gamma=0$, $\omega(0)=\pi^2/4$ and then the leading term in \eqref{215} becomes a special case of
\eqref{stst}. Also, if $\gamma=2$, $\omega(2)=0$, with $\omega>0$ for $\gamma<2$ and $\omega<0$ for $\gamma>2$.

From \eqref{3rdr} we see that $|\tm|$ is exponentially small, which implies exponentially large time scales. While the result \eqref{3rdr}
 is established analytically in Section \ref{sec4}, it can also be seen as a consequence of the following result. Let $\Rightarrow$ denote convergence in distribution.

 \begin{proposition}[Exponential limit law]\label{propren}
 Let $V$ be an exponential random variable with unit mean. Then,
 \eqan{\label{erg}
C \e^{-\beta b}\pas{b} \Rightarrow V, \quad {\rm as} \ b\ra\infty
}
with
 \eqan{\label{CCC}
C =\left(\frac{1}{\beta^2} + \frac{1}{\beta}\int_{0}^{\infty}e^{\beta u - u^2/2}du\right)^{-1}
}
 \end{proposition}
As mentioned in \cite{maglaraszeevi}, Proposition \ref{propren} can be established using 
a limit theorem from the theory of regenerative processes. In Section \ref{secnew} we give two proofs of Proposition \ref{propren}. The first proof is probabilistic and uses the theory of regenerative processes as pointed out in  \cite{maglaraszeevi}, and the second proof is analytic and uses the exact expressions for the Laplace transform.

We next give some results for the double limits $b\ra \infty$ with $\beta\ra\pm\infty,$ and also $b\ra 0$ with $\beta\ra -\infty$.
\begin{proposition}[Regime (ii)]\label{prop5a}
If $\beta\ra -\infty$ and $b\rightarrow \infty$, then
\eqan{\label{a1}
\tm \sim -\frac{1}{4}\beta^2-\frac{\pi^2}{b^2}\left[1+\frac{2}{b}\left(\frac{2}{-\beta}\right)^{1/3}\frac{{\rm Ai}(0)}{{\rm Ai}'(0)}\right],
}
where ${\rm Ai}(\cdot)$ is the Airy function.
If $\beta\ra \infty$ and $b\rightarrow +\infty$, then
\eqan{\label{a2}
\tm \sim -\frac{\beta}{\sqrt{2\pi}}\e^{-b\beta}\e^{-\beta^2/2}.
}
\end{proposition}

\begin{proposition}[Regime (iii)]\label{prop5b}
If $b\rightarrow 0$ with $\beta\ra-\infty$, we let $\Bs=b(-\beta)^{1/3}$ and for $\Bs$ fixed,
\eqan{\label{a3}
\tm \sim -\frac{1}{4}\beta^2-\left(-\frac{\beta}{2}\right)^{2/3}\eta,
}
where $\eta=\eta(\Bs)$ is the minimal solution to
\eqan{\label{a4}
\sqrt{-\eta}\cot\left[\Bs 2^{-1/3}\sqrt{-\eta}\right]=\frac{{\rm Ai}'(\eta)}{{\rm Ai}(\eta)}.
}
\end{proposition}
We comment that \eqref{a1} remains valid for $\beta\ra-\infty$ with $b>0$ fixed, \eqref{a2} remains valid for $\beta\ra\infty$ with $b>0$ fixed and can be obtained as a limiting case of \eqref{3rdr} (for $\beta\ra\infty$).

When $\Bs\ra\infty$ it follows from \eqref{a4} that $\eta\ra 0$ with $\eta\sim -2^{2/3}\pi^2 \Bs^{-2}$, while if $\Bs\ra 0^+$ we have $\eta\sim r_0=-2.338\ldots$ where $r_0$ is the least negative root of ${\rm Ai}(z)=0$. If $\eta=r_*<0$ where $r_*$ is the least negative root of ${\rm Ai}'(z)=0$, then $\Bs=\frac{1}{2}\pi 2^{1/3}(-r_*)^{-1/2}$.

\section{Proofs of Theorem \ref{thmneg} and Theorem \ref{thmpos}}\label{secc3}
From \eqref{diffusioneq} and \eqref{lst} it follows that the Laplace transform $\lst{x,b}(\theta)=Q(x;\theta)$ satisfies the ODE
\eqan{\label{31}
Q_{xx}+A(x)Q_x=\theta Q, \quad x<b,
}
with the boundary condition $Q(b;\theta)=1$ and the interface conditions $Q(0^+;\theta)=Q(0^-;\theta)$ and
$Q_x(0^+;\theta)=Q_x(0^-;\theta)$. For $x>0$ we have $A(x)=-\beta$ and then \eqref{31} admits solutions of the form $\e^{\alpha x}$ where  $\alpha=\alpha_{\pm}=\frac{1}{2}[\beta\pm\sqrt{\beta^2+4\theta}]$. For $x<0$, $A(x)=-x-\beta$ and then \eqref{31} becomes the Hermite equation, and the only solution that decays as $x\ra -\infty$ is proportional to $\e^{x^2/4}\e^{\beta x/2}D_{-\theta}(-x-\beta)$. It follows that
\eqan{\label{32}
Q(x;\theta)=k_0\exp\left(\frac{x^2}{4}+\frac{\beta x}{2}\right)\frac{D_{-\theta}(-\beta-x)}{D_{-\theta}(-\beta)}, \quad x<0,
}
and
\eqan{\label{33}
Q(x;\theta)=k_1{\rm e}^{\alpha_+(\theta)(x-b)}+(1-k_1){\rm e}^{\alpha_-(\theta)(x-b)}, \quad 0<x<b,
}
where $k_0$ and $k_1$ are independent of $x$. Here we wrote the solution in \eqref{33} in such a way so to automatically satisfy the boundary condition $Q(b;\theta)=1$. To determine $k_0$ and $k_1$ we can use the interface conditions at $x=0$, which imply that
\eqan{\label{34}
k_0=k_1\e^{-\alpha_+(\theta) b}+(1-k_1)\e^{-\alpha_-(\theta) b}
}
and
\eqan{\label{35}
-k_0\frac{D_{-\theta}'(-\beta)}{D_{-\theta}(-\beta)}=\frac{1}{2}\sqrt{\beta^2+4\theta}\left[k_1\e^{-\alpha_+(\theta) b}-(1-k_1)\e^{-\alpha_-(\theta) b}\right].
}
Solving the algebraic system in \eqref{34} and \eqref{35} for $k_0$ and $k_1$ leads to the expressions in Theorems \ref{thmneg} and \ref{thmpos}.

\section{Brief derivation of Propositions \ref{prop4}, \ref{prop5a} and \ref{prop5b}}\label{sec4}
We discuss the various asymptotic formulas for $\tm$. The expressions follow from routine manipulations of the parabolic cylinder functions that appear in \eqref{M} and \eqref{key}. Consider the solution to $D_{-\theta}(-\beta)M(\theta;\beta,b)=0$ which is equivalent to
\eqan{\label{41}
D_{-\theta}(-\beta)\cosh\left(\frac{b}{2}\sqrt{\beta^2+4\theta}\right)=
\frac{2}{\sqrt{\beta^2+4\theta}}D_{-\theta}'(-\beta)\sinh\left(\frac{b}{2}\sqrt{\beta^2+4\theta}\right).
}
We analyze this transcendental equation in various limits and find the least negative root $\tm$, only sketching the main points in the calculations.

For $b\ra\infty$ and $\beta>0$ we use the Taylor expansion, for $\theta\ra 0$,
\eqan{\label{}
\frac{D_{-\theta}'(-\beta)}{D_{-\theta}(-\beta)}=\frac{\beta}{2}-\theta\e^{\beta^2/2}\int_{-\infty}^\beta \e^{-u^2/2}{\rm d}u+O(\theta^2)
}
and, for $\theta\ra 0$ and $b\ra \infty$,
\eqan{\label{4333}
\coth(\frac{b}{2}\sqrt{\beta^2+4\theta})=1+2\e^{- b\beta}[1+O(\theta,\e^{-b\beta})].
}
Also, $\sqrt{\beta^2+4\theta}=(2/\beta)[1-2\theta/\beta^2+O(\theta^2)]$, so that \eqref{41} is equivalent to
\eqan{
\left[\frac{\beta}{2}-\theta\e^{\beta^2/2}\int_{-\infty}^\beta \e^{-u^2/2}{\rm d}u+O(\theta^2)
\right]
\frac{2}{\beta}\left[1-\frac{2\theta}{\beta^2}+O(\theta^2)\right]
=1+2 \e^{-b\beta} +O(\theta,\e^{-2b\beta})
}
and hence
\eqan{
-2\theta\left[\frac{1}{\beta^2}+\frac{1}{\beta}\e^{\beta^2/2}\int_{-\infty}^\beta \e^{-u^2/2}{\rm d}u\right]\sim 2\e^{-b\beta},
}
which leads to the exponentially small $\tm$ in \eqref{3rdr}.
A similar calculation applies for fixed $b>0$ and $\beta\ra\infty$, which leads to \eqref{a2} in Proposition \ref{prop5a}.

For $b\ra\infty$ with $\beta<0$ the solution $\tm$ will be close to the apparent branch point at $\theta=-\beta^2/4$. From \eqref{41} if we set $\theta=-\beta^2/4-\Omega/b^2$ we obtain
\eqan{\label{42}
\frac{D_{-\theta}'(-\beta)}{D_{-\theta}(-\beta)}+O(\theta+\beta^2/4)=\frac{b}{\sqrt{\Omega}}\tan(\sqrt{\Omega}).
}
If $b\ra\infty$ with a fixed $\beta<0$, $\sqrt{\Omega}$ must be close to a zero of the tangent function, so that $\Omega\sim\pi^2/4$ (for the least negative solution $\tm$). Then estimating the difference $\sqrt{\Omega}-\pi/2$ using \eqref{42} leads to \eqref{214}.

If $b\ra\infty$ with $\beta\ra 0$ we again obtain \eqref{42} and for small $\beta$ we can further approximate
\eqan{\label{}
\frac{D_{-\theta}'(-\beta)}{D_{-\theta}(-\beta)}\sim\frac{2}{\beta}=\frac{2b}{\gamma}.
}
Then \eqref{stst} follows from \eqref{42} with $\Omega$ replaced by $\omega=b^2\tm -\gamma^2/4$. If $\beta=0$ we can express the parabolic cylinder functions $D_{-\theta}(0)$ and $D_{-\theta}'(0)$ in terms of Gamma functions, and \eqref{41} becomes
\eqan{\label{43}
\frac{\tanh(b\sqrt{\theta})}{\sqrt{\theta}}=\frac{-\Gamma(\frac{\theta}{2})}{\sqrt{2}\Gamma(\frac{\theta+1}{2})}.
 }
 Then \eqref{215} follows by solving \eqref{43} for $\tm$ with $b\ra\infty$, where $\tm\sim-\pi^2/(4b^2)=O(b^{-2})$. Then $\Gamma(\frac{\theta+1}{2})\sim\sqrt{\pi}$ and $\Gamma(\frac{\theta}{2})\sim\frac{2}{\theta}$.

If $\beta\ra-\infty$ and $b\ra\infty$ or $b=o(1)$, \eqref{a1} follows by a calculation similar to \eqref{42}, except now $D_{\beta^2/4}(-\beta)$ becomes proportional to ${\rm Ai}(0)$ in this limit. Finally, to obtain \eqref{a3} we use the approximation
\eqan{\label{44}
\frac{-\frac{{\rm d}}{{\rm d}\beta}D_{-\theta}'(-\beta)}{D_{-\theta}(-\beta)}\sim \left(\frac{-\beta}{2}\right)^{1/3}\frac{{\rm Ai}'(\eta)}{{\rm Ai}(\eta)}, \quad \eta=\frac{\theta+\beta^2/4}{(-\beta/2)^{2/3}}
}
which applies for $-\beta\ra\infty$ and $-\theta=\beta^2/4+O(|\beta|^{2/3})$. With \eqref{44}, \eqref{a4} follows from \eqref{41}. This concludes the sketched derivation of Propositions \ref{prop4}, \ref{prop5a} and \ref{prop5b}.

\section{Two proofs of Proposition \ref{propren}}\label{secnew}
\subsection{Probabilistic proof}
The crucial observation is that the one-dimensional Halfin-Whitt diffusion process is an ergodic one-dimensional diffusion process that has, by the strong Markov property, the origin as a regeneration point. To make this formal, let $y$ be a fixed positive number, where $y < b$.  This value can be used to construct a sequence of random times $\{T_{k}\}_{k \geq 1}$, where
\begin{eqnarray} T_1 = \inf\{t \geq 0: X(t) = 0, \sup_{0 \leq s \leq t}X(s) \geq y\} \end{eqnarray} and for each $n \geq 1$,
\begin{eqnarray} T_{n+1} = \inf\{t \geq T_n: X(t) = 0, \sup_{T_n \leq s \leq t}X(s) \geq y\}. \end{eqnarray}  The Halfin-Whitt diffusion is a regenerative process with respect to these regeneration times, which form a delayed renewal process that, with probability one, has a finite number of points in each compact interval. To make matters simpler to state, we assume that $X(0) = 0$, but the procedure outlined here can be adjusted for any arbitrary initial condition.

Define $a(b)$ to be the probability that our process gets above level $b$ in the random interval $[0, T_1]$.  Due to the Strong Markov property, along with the fact that $y < b$, we observe that this probability is just the probability that a Brownian motion with drift $-\beta$ and diffusion coefficient $2$ reaches level $b$ before level $0$.  Hence,
\begin{eqnarray} a(b) = \frac{e^{\beta y} - 1}{e^{\beta b} - 1}.
\end{eqnarray}
Moreover, the expected length of each regenerative cycle is just the expected amount of time it takes the Halfin-Whitt diffusion to reach level $y$, starting from $0$, plus the expected amount of time it takes the diffusion to go from level $y$ back to $0$.  In other words,
\begin{eqnarray} \expect{T_1} = (e^{\beta y} - 1)\left[\frac{1}{\beta^2} + \frac{1}{\beta}\int_{0}^{\infty}e^{\beta u - u^2/2}du\right] \end{eqnarray}  Therefore, by \cite[Theorem 4.2 on p. 181]{asmussen}, we see that as $b \to \infty$,
\begin{eqnarray} \frac{a(b)}{\expect{T_1}}T_{x}(b) \Rightarrow V \end{eqnarray} where $V$ is an exponential random variable with rate one.  Moreover, observe that for each $y < b$
\begin{eqnarray}
\frac{a(b)}{\expect{T_1}} &=& \frac{1}{(e^{\beta b} - 1)\left[\frac{1}{\beta^2} + \frac{1}{\beta}\int_{0}^{\infty}e^{\beta u - u^2/2}du\right]} \end{eqnarray}
which does not depend on $y$.  Hence,  as $b \to \infty$, we obtain \eqref{erg}.
A very similar asymptotic result carries through as well for $\expect{T_{x}(b)}$, due to Asmussen \cite[Proposition 4.1 on p.~180]{asmussen}.
From \eqref{erg} it follows that the time it takes to reach a level $b$ is roughly exponential in $b$, which says that extreme congestion is not observed on relatively short time scales. For $\beta>0$ the mean first passage time is
$
\expect{\pas{b}}\sim 1/|\tm|
$ for $x$ bounded away from $b$ and $b\ra\infty$.

\subsection{Analytic proof}
Consider the Laplace transforms in \eqref{28} and \eqref{29} on the scale $\theta=O(|\theta_{\rm max}|)=O(\e^{-\beta b})$ in \eqref{3rdr}, and then scale time as $t=T/|\theta_{\rm max}|=O(\e^{b\beta})$. Since $D_0(-\beta)=\e^{-\beta^2/4}$ we have
\eqan{\label{52a}
\frac{D_{-\theta}'(-\beta-x)}{D_{-\theta}(-\beta)}\sim \exp\left(-\frac{x\beta}{2}-\frac{x^2}{4}\right), \quad \theta\to 0.
}
Then we write $M$ in \eqref{M} as
 \eqan{\label{52b}
M=\frac{\sinh\left(\frac{b}{2}\sqrt{\beta^2+4\theta}\right)}{\sqrt{\beta^2+4\theta}}\Big[\sqrt{\beta^2+4\theta}\coth\Big(\frac{b}{2}\sqrt{\beta^2+4\theta}\Big)-
\frac{2D_{-\theta}'(-\beta)}{D_{-\theta}(-\beta)}\Big].
}
For $\theta\to 0$ we use \eqref{4333}, and the recurrence relation for the parabolic cylinder function $D_p'(z)=-\frac z 2 D_p(z)+p D_{p-1}(z)$, which for $p\to 0$ yields
\eqan{\label{52c}
\frac{D_{p}'(z)}{D_{p}(z)}=-\frac z 2 +p\frac{D_{-1}(z)}{D_0(z)}+O(p^2).
}
Then \eqref{52b} becomes
 \eqan{\label{52d}
M&=\frac{1}{2\beta}\e^{b\beta/2}\Big[(1+2\e^{-b\beta}+O(\theta\e^{-b\beta},\e^{-2 b\beta}))(\beta+\frac{2\theta}{\beta}+O(\theta^2))
-\beta+2\theta\e^{\beta^2/4} D_{-1}(-\beta)+O(\theta^2)\Big]\nn\\
&\sim \frac{1}{\beta}\e^{b\beta/2}\Big(\theta\Big[\frac{1}{\beta}+\e^{\beta^2/4}D_{-1}(-\beta)\Big]+2\beta\e^{-b\beta}\Big)\nn\\
&=\e^{-b\beta/2}\Big(\frac{\theta+|\theta_{\rm max}|}{|\theta_{\rm max}|}\Big),
}
$\theta=O(|\theta_{\rm max}|)$, as $D_{-1}(-\beta)=\e^{\beta^2/4}\int_{-\infty}^\beta \e^{-u^2/2}d u= \e^{-\beta^2/4}\int_{0}^\infty \e^{\beta u} \e^{-u^2/2}d u$. Then inverting the transform in \eqref{28}, using \eqref{52a} and \eqref{52d}, yields (with ${\rm Br}$ a vertical Bromwich contour in the half-plane $\Re(\xi)>0$)
\eqan{\label{52e}
P(x,t)\sim\frac{|\theta_{\rm max}|}{2\pi i} \int_{\rm Br}\frac{1}{1+\xi}\e^{\xi T}d\xi=|\theta_{\rm max}|\e^{-T},
}
since $\theta t = \xi T$ if $\theta=|\theta_{\rm max}| \xi$.
This yields the exponential limit law for $x<0$.

If $x>0$ so that both $b\to\infty$ and $b-x\to\infty$, then an analogous calculation using \eqref{29} again leads to \eqref{52e}. However, if $x,b\to\infty$ in such a way that $b-x=O(1)$, then \eqref{29} leads to the limit
\eqan{\label{}
\lst{x,b}(\theta)\to \e^{(x-b)\beta}+ \frac{|\theta_{\rm max}|}{\theta+|\theta_{\rm max}|}\Big[1-\e^{(x-b)\beta}\Big]
}
and this inverts to, on time scales $t=O(\e^{b\beta})$,
\eqan{\label{52f}
P(x,t)&\sim\frac{|\theta_{\rm max}|}{2\pi i} \int_{\rm Br}\Big[\e^{(x-b)\beta}+\frac{1-\e^{(x-b)\beta}}{1+\xi}\Big]e^{T\xi}d\xi\nn\\
&=|\theta_{\rm max}|\Big[\e^{(x-b)\beta}\delta(T)+(1-\e^{(x-b)\beta})\e^{-T}\Big].
}
The term proportional to $\delta(T)$ represents a probability mass on the large time scale, which corresponds to sample paths that hit $b$ in a short time (at least $t=o(e^{b\beta})$). The actual density of $P(x,t)$ does not have mass at $t=0$, and finding the expansion of the density on shorter time ranges would require a different asymptotic analysis, and a different approximation to $M$. We refer the reader to \cite{A,B}, where such problems were analyzed in detail for some simpler models. There both exponentially large and $t=O(1)$ time scales were considered, and they were related to one another by asymptotic matching.

\section*{Acknowledgments}
The work of Charles Knessl was supported partially by NSA grants H 98230-08-1-0102 and H 98230-11-1-0184. The work of Johan van Leeuwaarden was supported by an ERC Starting Grant.
\begin{footnotesize}

\end{footnotesize}
\end{document}